\documentclass[11pt]{amsart}
\usepackage{amsmath, amssymb, amsthm, amsfonts,enumitem,tikz,relsize,array}
\usetikzlibrary{math}
\usepackage[normalem]{ulem}
\usepackage[colorlinks = true, citecolor = blue, hypertexnames=false]{hyperref}
\usepackage{xcolor}
\hypersetup{
    colorlinks,
    linkcolor={red!70!black},
    citecolor={blue!70!black},
    urlcolor={blue!80!black}
}
\usepackage{framed,booktabs,float,caption,scalerel}
\usepackage[OT1]{fontenc}
\usepackage{comment}
\usepackage{diagbox}
\usepackage{multicol}
\usepackage[textsize=scriptsize, bordercolor=lightgray, linecolor=lightgray, backgroundcolor=white]{todonotes}

\theoremstyle{plain}
\newcounter{MainTheoremCounter}

\newtheorem{MainTheorem}[MainTheoremCounter]{Theorem}
\newtheorem{MainCorollary}{Corollary}[MainTheoremCounter]

\newtheorem{theorem}{Theorem}[section]
\newtheorem{lemma}[theorem]{Lemma}
\newtheorem{proposition}[theorem]{Proposition}
\newtheorem{fact}[theorem]{Fact}
\newtheorem*{theorem*}{Theorem}
\newtheorem*{lemma*}{Lemma}
\newtheorem*{corollary*}{Corollary}

\newcommand{\thistheoremname}{}
\newtheorem*{genericthm}{\thistheoremname}
\newenvironment{namedthm}[1]
  {\renewcommand{\thistheoremname}{#1}%
   \begin{genericthm}}
  {\end{genericthm}}  

\theoremstyle{definition}
\newtheorem{example}[theorem]{Example}
\newtheorem{remark}[theorem]{Remark}
\newtheorem*{remark*}{Remark}

\newtheorem{question}{Question}
\newtheorem*{definition*}{Definition}
\newtheorem*{setup}{Setup}
\newtheorem*{notation}{Notation}

\theoremstyle{remark}
\newtheorem{claim}{Claim}

\newenvironment{proofclaim}[1][Proof of Claim]{\begin{proof}[#1]}{\end{proof}}

\newcommand{\qoplus}{\mathbin{\scalerel*{\left(\mkern-2mu+\mkern-2mu\right)}{\bigcirc}}}

\DeclareMathOperator{\ad}{ad}
\DeclareMathOperator{\tr}{tr}
\DeclareMathOperator{\im}{im}

\DeclareMathOperator{\charac}{char}
\DeclareMathOperator{\Aut}{Aut}
\DeclareMathOperator{\End}{End} 
\DeclareMathOperator{\Sym}{Sym} 
\DeclareMathOperator{\Alt}{Alt}
\DeclareMathOperator{\icosa}{Icosa}
\DeclareMathOperator{\SO}{SO}

\newcommand{\abs}[1]{\lvert#1\rvert}
\newcommand{\N}{\mathbb{N}}
\newcommand{\Z}{\mathbb{Z}}
\newcommand{\R}{\mathbb{R}}

\newcommand{\cU}{\mathcal{U}}
\DeclareMathOperator{\Ar}{Ar}
\DeclareMathOperator{\Ob}{Ob}

\begin{document}
\title[Abstract 3D-rotation groups]{Abstract 3D-rotation groups and recognition of icosahedral modules}
\author{Lauren  McEnerney}
\address{Department of Mathematics and Statistics\\
California State University, Sacramento\\
Sacramento, CA 95819, USA}
\email{lmcenerney@csus.edu}
\author{Joshua Wiscons}
\address{Department of Mathematics and Statistics\\
California State University, Sacramento\\
Sacramento, CA 95819, USA}
\email{joshua.wiscons@csus.edu}
\date{\today}

\begin{abstract}
We introduce an abstract notion of a 3D-rotation module for a group $G$ that does not require the module to carry a vector space structure, a priori nor a posteriori. We prove that, under an expected irreducibility-like assumption, the only finite $G$ with such a module are those already known from the classical setting: $\Alt(4)$, $\Sym(4)$, and $\Alt(5)$. Our main result then studies the  module structure when $G =  \Alt(5)$ and shows that, under certain natural restrictions, it is fully determined and generalizes that of the classical icosahedral module.

We  include an application to the recently introduced setting of modules 
with an additive dimension, a general setting allowing for simultaneous treatment of classical representation theory of finite groups as well as representations  within various well-behaved model-theoretic settings such as the $o$-minimal and finite Morley rank ones.
Leveraging our recognition result for icosahedral modules, we classify the faithful $\Alt(5)$-modules with additive dimension that are dim-connected of dimension $3$ and without $2$-torsion.
\end{abstract}

\maketitle

\section{Introduction}

The classification of the finite subgroups of $\SO_3(\R)$ is a  classic of mathematics with broad applications throughout and beyond the discipline. Known in spirit to Plato (and at least in part well before), modern approaches to its proof can be traced back to Klein~\cite{KlF93} with many variations through  different lenses appearing over the years. 

Here we revisit the classification in a significantly more general context  that allows for actions on groups that need not be vector spaces. 
Our initial goal  was to capture modules that---while perhaps not finite dimensional over a field---live in a universe with a well-behaved notion of dimension, such as the model-theoretic settings of finite Morley rank or $o$-minimality, and it was a desire to classify minimal modules in these dimensional contexts that actually started this project. The paper~\cite{CDWY24} (echoing~\cite{CDW23}) asked for a classification of the minimal faithful $\Alt(5)$-modules with an additive dimension in characteristics other than $2$ and $5$, and this essentially appears  as our Corollary~\ref{cor.Alt5Recognition-dim}, which also includes characteristic $5$.
But as the project evolved, we realized that our main result, Theorem~\ref{thm.Alt5Recognition}, could be obtained without  assuming the existence of a dimension.

We now define our setting. Regarding notation: for $G$ a group acting on a group $V$, we use centralizer and normalizer notation consistent with the semidirect product $G\ltimes V$. Actions will  be on the left.

\begin{definition*}
Let $G$ be a group. A $G$-module $V$ will be called a \textbf{3D-rotation module} if every nontrivial $g\in G$ can be assigned a nontrivial subgroup $A_g \le V$, called the \textbf{axis} of $g$, such that the following hold for all $h\in G\setminus \{1\}$:
\begin{enumerate}
    \item $hA_g = A_{hgh^{-1}}$;
    \item\label{item.def.3Drotation-C} $h\in C_G(A_g)$ if and only if $A_h = A_g$;
    \item\label{item.def.3Drotation-N} $h\in N_G(A_g)\setminus C_G(A_g)$ implies $h$ inverts $A_g$.
\end{enumerate}
In the case that some axis is $G$-invariant, the module is called \textbf{axial}.
\end{definition*}

To  capture the classical setting when $G\le \SO_3(\mathbb{R})$, one may choose the axis $A_g$ to be as large as possible: the full centralizer in $V = \mathbb{R}^3$ of $g$, i.e.~$A_g = C_V(g)$. In fact, in this case, the axis can equivalently be defined as the image of the \emph{trace map} $\tr_g = 1 + g + \cdots + g^{n-1}$ where $n = |g|$; this will be our preferred point of view. The trace map is reviewed in \S\S~\ref{s.trace-adjoint}. 

Although our definition does not explicitly mention dimension,  items~\eqref{item.def.3Drotation-C} and \eqref{item.def.3Drotation-N} encode that axes behave as though they are $1$-dimensional. Also, that group elements act as rotations is only addressed for involutions  (in~\eqref{item.def.3Drotation-N}), but this turns out to be enough to effectively explore the nonaxial case. 

\begin{remark*}
The external property of having a 3D-rotation module imposes some quite familiar internal restrictions on $G$: taking $C = C_G(A)$ for $A$ any axis, one finds that $C$ is a \textsc{ti}-subgroup (i.e.~the intersection of $C$ with any distinct conjugate is trivial), and moreover, $C$ is of index at most $2$ in its normalizer. In total, when $C$ is proper in $G$, we have a quasi-Frobenius---and perhaps  Frobenius---pair $C <  G$ in the sense of~\cite{ZaS23} (see also ~\cite{ACD24}).
\end{remark*}
\subsection{Results}\label{s.results}
Our first theorem is the expected analog of the  classification of the finite subgroups  of $\SO_3(\R)$. Its statement and proof are the natural adaptations of Klein's work to our more general context.

\begin{MainTheorem}[Generalizing Klein]\label{thm.ClassificationRotationGroups}
Let $G$ be a nontrivial finite group and $V$ a nonaxial 3D-rotation module for $G$. Then $G$ is isomorphic to one of $\Alt(4), \Sym(4), \Alt(5)$ and the centralizer of each axis is cyclic. 
\end{MainTheorem}


Theorem~\ref{thm.ClassificationRotationGroups} leaves  open the obvious question of the actual $G$-module structure on $V$. Our second---and main---theorem  addresses this for $\Alt(5)$ showing that, under certain natural restrictions, it must be icosahedral.

\begin{definition*}[Icosahedral Modules]
Let $R = \Z[\phi]$ with $\phi$ satisfying $\phi^2 = \phi +1$. Fix a generating set $\{\alpha,\gamma\}$ for  $\Alt(5)$ with $|\alpha| = 2$,  $|\gamma| = 3$, and $|\alpha\gamma| = 5$. Define $\icosa(R)$ (with respect to $\alpha,\gamma$) to be the $R[\Alt(5)]$-module $Re_1 \oplus Re_2 \oplus Re_3$ with action given as follows:
\begin{align*}
\alpha e_1 & = -e_1 + \phi e_3 & \gamma e_1 &= e_3\\
\alpha e_2 & = -e_2 + \phi e_3 & \gamma e_2 &= e_1\\
\alpha e_3 & = e_3 & \gamma e_3 &= e_2
\end{align*}
Then, for any  $R$-module $L$ without $2$-torsion, set $\icosa(L) = \icosa(R) \otimes_R L$ with $\Alt(5)$ acting trivially on $L$.
\end{definition*}

Our definition for $\icosa(R)$ is taken from the Atlas \cite{atlas3}. We find this representation easy to work with, but do note that the corresponding matrix for $\alpha$ is \emph{not} in $\SO_3(R)$.

\begin{remark*}
We insist  $L$ have no $2$-torsion to rule out  modules $V = \icosa(L)$ for which $C_V(\Alt(5))\neq 0$. Indeed, if $\ell\in L$ is an involution, then $(e_1+e_2+e_3)\otimes \ell$ is centralized by $\alpha$ and $\gamma$, hence by $\Alt(5)$. Conversely, if $w = \sum e_i\otimes \ell_i$ is in $C_V(\Alt(5))$, then as $w$ is fixed by $\gamma$, $\ell_1=\ell_2=\ell_3$, and as $w$ is also fixed by $\alpha$,  $e_1\otimes 2\ell_1 + e_2\otimes 2\ell_1 + e_3\otimes (1-2\phi)\ell_1 = 0$, implying  $2\ell_1 =0$.
\end{remark*}

The class of all $\icosa(L)$ captures a number of modules that differ in varying degrees from the classical setting. In particular, one could take $L$ to be a sum of fields of different characteristics, or, more interestingly, take $L$ to be a Pr\"ufer $p$-group as in the following example.  

\begin{example}
For $p$ a prime different than $2$ or $5$, the Pr\"ufer $p$-group $C_{p^\infty}$ has an automorphism $\phi$ satisfying $\phi^2 = \phi +1$ precisely when $p\equiv \pm1$ (mod $5$). Indeed,  $\Aut(C_{p^\infty})$ is isomorphic to the $p$-adic integers, so determining when such $\phi$ exists reduces, by Hensel's lemma, to determining when $x^2 - x - 1$ has root modulo $p$. This in turn reduces to determining when $5$ is a square modulo $p$, and quadratic reciprocity then yields the desired condition. Thus, we have modules of the form $\icosa(C_{p^\infty})$ whenever $p\equiv \pm1$ (mod $5$).
\end{example}

Theorem~\ref{thm.Alt5Recognition} is a recognition theorem for a certain class of icosahedral modules, with the main hypothesis being a local and internal version of the 3D-rotation module axioms with respect to a specific choice of axes. Our approach is direct and constructive, and it appears to be new even in the classical setting. 

The hypotheses of Theorem~\ref{thm.Alt5Recognition}  reference the characteristic of a module (if it exists) as well as the notion of $p$-divisibility, which we take to be a robust form of characteristic \emph{not} $p$. 
For $V$  an abelian group and $p$ a prime, we say $V$ has \textit{characteristic $p$} if $V$ is an elementary abelian $p$-group. Also, $V$ is \textit{$p$-divisible} if for all $v\in V$, there is $w\in V$ such that $v = pw$; if there is a unique such $w$ for each $v$, then $V$ is \textit{uniquely $p$-divisible}. If $V$ is $p$-divisible for all primes $p$, then $v = mw$ has a solution for all  $m\in \Z$, in which case $V$ is simply said to be \textit{divisible}.

\begin{MainTheorem}\label{thm.Alt5Recognition}
Let $V$ be a uniquely $2$-divisible $\Alt(5)$-module that either has prime characteristic or is divisible. 
Whenever $V$ is $|g|$-divisible for nontrivial $g\in \Alt(5)$,  assume that the elements of $N_{\Alt(5)}(\langle g \rangle) \setminus \langle g \rangle$ invert  $\tr_g(V)$.
Then 
$V\cong \icosa(L)$ with the $\Z[\phi]$-module $L$ defined as follows: 
\begin{itemize}
\item choose any involution $\alpha\in \Alt(5)$, and set $L = \tr_\alpha(V)$;
\item choose either conjugacy class $\mathcal{C}$ of $5$-cycles, and let 
$\phi$ be the image in $\End(L)$ of $\frac{1}{4}\sum_{g\in \mathcal{C}} g$.
\end{itemize}
\end{MainTheorem}

With Theorem~\ref{thm.ClassificationRotationGroups} in hand, the ``normalizer assumption'' in the statement of Theorem~\ref{thm.Alt5Recognition} is easily seen to hold of nonaxial 3D-rotation modules for $\Alt(5)$ with respect to the choice  $A_g = \tr_g(V)$, so this  yields the classification of such modules assuming they have prime characteristic or are divisible. 
The additional hypothesis that $V$ be uniquely $2$-divisible is more-or-less expected since 3D-rotation modules of characteristic $2$ are axial: 
 if $\charac V = 2$, then $N_G(A_g) = C_G(A_g)$ for all $g$, which Theorem~\ref{thm.ClassificationRotationGroups} implies is an axial configuration.

\begin{MainCorollary}\label{cor.Alt5Recognition-3Drot}
Let $V$ be a nonaxial 3D-rotation module for $\Alt(5)$ with respect to the axes $A_g = \tr_g(V)$. If $V$ is uniquely $2$-divisible and  either has prime characteristic or is divisible, then the conclusion of Theorem~\ref{thm.Alt5Recognition} holds.
\end{MainCorollary}

Theorem~\ref{thm.Alt5Recognition} also allows us to essentially complete the classification---started in~\cite{CDWY24}---of the minimal $\Alt(5)$-modules with an additive dimension. 
Background on this quite natural setting is provided in \S\S~\ref{s.Alt5Recognition-dim}.

\begin{MainCorollary}\label{cor.Alt5Recognition-dim}
Let $V$ be a faithful $\Alt(5)$-module with an additive dimension. Assume $V$ is dim-connected and $\dim V = 3$. If $V$ is without $2$-torsion, then  the conclusion of Theorem~\ref{thm.Alt5Recognition} holds and $\dim L = 1$.
\end{MainCorollary}

This result emphasizes---like those in \cite{CDW23,CDWY24}---that although the setting of modules with an additive dimension allows  for a variety of non-classical structures, the minimal ones contain no real surprises. If one specializes to the finite Morley rank setting, then this point becomes particularly sharp in light of recent work by Borovik~\cite{BoA24}.

\subsection{Questions}

We  record a handful of questions, many of which are implicit above. We have not thought deeply on them so make no claim on their difficulty. With luck at least one will land properly between easy and hopeless. All can be taken as is or specialized to a setting where the module carries a reasonable notion of dimension. 

\begin{question}
Is there a version of Corollary~\ref{cor.Alt5Recognition-3Drot} for $\Alt(4)$ and $\Sym(4)$? 
\end{question}

\begin{question}
What---if anything---can be said about finite groups with a faithful \emph{axial} 3D-rotation module? Must they be solvable? 
\end{question}

\begin{question}
What can be said about \emph{infinite} groups with a nonaxial 3D-rotation module? What are  interesting examples that do not already occur in some $\SO(V)$ with $V$ a $3$-dimensional vector space? 
\end{question}

\begin{question}
Is there a  generalization of  abstract 3D-rotation modules to higher dimensions for which a version of Theorem~\ref{thm.ClassificationRotationGroups} holds? 
\end{question}

\section{Classification of the finite 3D-rotation groups}

We work towards a proof of Theorem~\ref{thm.ClassificationRotationGroups}. 
Our approach very closely follows Breuillard's rendition~\cite{BrE23} of the classic work of Klein~\cite{KlF93}.

Notice that if $V$ is a 3D-rotation module for $G$ then $C_G(A)$ has index at most two in $N_G(A)$ for each axis $A$. Thus, axes come in two flavors: we say $A$ has \textit{positive type} if $N_G(A) = C_G(A)$ and \textit{negative type} otherwise. 

We adopt the following setup for the remainder of the section.

\begin{setup}
Let $G$ be a nontrivial finite group---with $G^*$ its nontrivial elements---and let $V$ be a nonaxial 3D-rotation module for $G$. Set $\mathbb{A}= \{A_g \mid g\in G^*\}$ with $\mathbb{A}^+$ and $\mathbb{A}^-$ denoting the subsets of positive- and negative-type axes, respectively. 
Choose representatives $A_1,\ldots A_r$ from each orbit of $G$ on $\mathbb{A}^+$ and similarly choose $B_1,\ldots B_s$ from the orbits of $G$ on $\mathbb{A}^-$.
\end{setup}

\subsection{Parameters on orbits}

Our analysis begins by pinning down the possibilities for the parameters $r$ and $s$. Classically, $(r,s)$ is one of $(1,0)$, $(1,1)$, or $(0,3)$ with $(1,0)$ corresponding to cyclic groups, $(1,1)$ to $\Alt(4)$ as well as dihedral groups of order $2n$ with $n$ odd, and $(0,3)$ to $\Sym(4)$, $\Alt(5)$, and dihedral groups of order $2n$ with $n$ even. As we assume the module is nonaxial from the outset, our list will be limited to $(1,1)$ and $(0,3)$.

\begin{lemma}[{cf.~\cite[Equation~(5-3)]{BrE23}}]\label{lem.OrbitEquation}
We have that \[ 1 - r - \frac{s}{2} = \frac{1}{|G|} - \sum_{i=1}^r \frac{1}{\abs{N_G(A_i)}} - \sum_{i=1}^s \frac{1}{\abs{N_G(B_i)}}.\]
\end{lemma}
\begin{proof}
The definition of a 3D-rotation module implies that $G^*$ is partitioned by the collection $\{C^*_G(A) \mid A\in \mathbb{A}\}$ (using  notation $H^* = H\setminus \{1\}$), so 
\[\abs{G} = 1 + \sum_{A\in \mathbb{A}} \abs{C^*_G(A)}.\] 
Working orbit-by-orbit and separating axes by type, we get
\[\abs{G} = 1 +  \sum_{i=1}^r \abs{G:N_G(A_i)}\cdot\abs{C^*_G(A_i)} +  \sum_{i=1}^s \abs{G:N_G(B_i)}\cdot\abs{C^*_G(B_i)},\]
which then yields
\[1 = \frac{1}{\abs{G}} + \sum_{i=1}^r \frac{\abs{C_G(A_i)} - 1}{\abs{N_G(A_i)}} + \sum_{i=1}^s \frac{\abs{C_G(B_i)} - 1}{\abs{N_G(B_i)}}.\]
The relationship between normalizers and centralizers of axes is completely determined by their type, so we obtain
\[1 = \frac{1}{\abs{G}} + \sum_{i=1}^r \left(1 - \frac{1}{\abs{N_G(A_i)}}\right) + \sum_{i=1}^s \left(\frac{1}{2} - \frac{1}{\abs{N_G(B_i)}}\right),\]
which implies the lemma.
\end{proof}

\begin{lemma}\label{lem.rsRestrict}
$(r,s)$ is equal to $(1,1)$ or $(0,3)$.
\end{lemma}

\begin{proof}
Note that if $A \in\mathbb{A}^+$, we have that $\abs{N_G(A)}=  |C_G(A)|\ge 2$, and when $A \in\mathbb{A}^-$, we get $\abs{N_G(A)}= 2\cdot |C_G(A)|\ge 4$. Thus, by Lemma~\ref{lem.OrbitEquation},
\[1 - r - \frac{s}{2} = \frac{1}{|G|} - \sum_{i=1}^r \frac{1}{\abs{N_G(A_i)}} - \sum_{i=1}^s \frac{1}{\abs{N_G(B_i)}} \ge \frac{1}{|G|} - \frac{r}{2} - \frac{s}{4}.\]
We conclude that $1 - \frac{r}{2} - \frac{s}{4} \ge \frac{1}{|G|} > 0$, which
insists $0 \leq r \leq 1$ and $0 \leq s \leq 3$. Moreover,  $(r,s)\neq(1,2),(1,3)$. 

Returning to Lemma~\ref{lem.OrbitEquation} and multiplying by $|G|$, we obtain 
\[ \left(1 - r - \frac{s}{2}\right)\abs{G} = 1 - \sum_{i=1}^r \dfrac{\abs{G}}{\abs{N_G(A_i)}} - \sum_{i=1}^s \dfrac{\abs{G}}{\abs{N_G(B_i)}},\]
where the right side is now a sum of integers. As $\mathbb{A}$ is nonempty, at least one of $r$ and $s$ is positive. 
Moreover, the module is nonaxial, so $N_G(A_i)$ and $N_G(B_i)$ are proper in $G$. Thus, $\left(1 - r - \frac{s}{2}\right)\abs{G} < 0$, implying that $r + \frac{s}{2} > 1$. Only $(r,s)=(1,1),(0,3)$ remain.
\end{proof}

\subsection{Proof of Theorem~\ref{thm.ClassificationRotationGroups}}

Our main tool for Theorem~\ref{thm.ClassificationRotationGroups} is Lemma~\ref{lem.OrbitEquation}, which becomes highly focused in light of Lemma~\ref{lem.rsRestrict}. 
The identification process is straightforward in the cases leading to the tetrahedral and octahedral groups as there is canonical orbit of axes to study. In the case leading to the icosahedral group, we instead consider the action on the Sylow $2$-subgroups, which serve as proxies for orthogonal triples of axes.

\begin{remark}\label{rem.BasicPrinciples}
The following trivial, but key, principle (already appearing in the proof of Lemma~\ref{lem.rsRestrict}) will be used repeatedly: if $A \in\mathbb{A}^+$, then $\abs{N_G(A)}=  |C_G(A)|\ge 2$, and if $A \in\mathbb{A}^-$, then $\abs{N_G(A)}= 2\cdot |C_G(A)|\ge 4$.

Also, as remarked earlier, the definition of a 3D-rotation module implies that for distinct $A_1,A_2\in \mathbb{A}$, $C_G(A_1) \cap C_G(A_2) = 1$. This in turn controls $N_G(A_1) \cap N_G(A_2)$ as it embeds into  $N_G(A_1)/C_G(A_1) \times N_G(A_2)/C_G(A_2)$, so $N_G(A_1) \cap N_G(A_2)$ has order dividing $4$.
\end{remark}

\begin{proof}[Proof of Theorem~\ref{thm.ClassificationRotationGroups}]
By Lemma~\ref{lem.rsRestrict}, we have two cases: $(r,s)=(1,1),(0,3)$.

\begin{claim}\label{cl.rs11}
If $(r,s) = (1,1)$, then $G\cong\Alt(4)$. 
\end{claim}
\begin{proofclaim}
    Choosing  $A\in \mathbb{A}^+$ and  $B\in \mathbb{A}^-$, Lemma~\ref{lem.OrbitEquation}  reads
    \[ \frac{1}{2} + \frac{1}{\abs{G}} = \frac{1}{\abs{N_G(A)}} + \frac{1}{\abs{N_G(B)}}, \]
    so as $\abs{N_G(B)} \geq 4$, $\frac{1}{2} < \frac{1}{\abs{N_G(A)}} + \frac{1}{4}$, implying that $2\le  \abs{N_G(A)} <4$.

    If $\abs{N_G(A)} = 2$, then $G = N_G(B)$, making $V$ axial, a contradiction. So $\abs{N_G(A)} = 3$, and then $\frac{1}{6} <   \frac{1}{\abs{N_G(B)}}$. Thus $4\le  \abs{N_G(B)} <6$, with $\abs{N_G(B)}$ even,  so $\abs{N_G(B)}=4$ and $\abs{G} = 12$. The group is now easily identified: the orbit containing $A$ has size $4$, and as  $N_G(A) = C_G(A)$, the action of $G$ on this orbit is readily seen to have a trivial kernel (see Remark~\ref{rem.BasicPrinciples}).
\end{proofclaim}

We may now assume $(r,s) = (0,3)$. Let $B_1, B_2, B_3$ be representatives of the three negative-type orbits ordered so that $\abs{N_G(B_1)}\le\abs{N_G(B_2)}\le \abs{N_G(B_3)}$. Lemma~\ref{lem.OrbitEquation} then gives
\[ \frac{1}{2} + \frac{1}{\abs{G}} = \frac{1}{\abs{N_G(B_1)}} + \frac{1}{\abs{N_G(B_2)}} + \frac{1}{\abs{N_G(B_3)}}.\]
Remember that $\abs{N_G(B_i)}$ is even. As $\abs{N_G(B_1)}\ge 4$, we find $\abs{N_G(B_1)} = 4$ and  $\abs{N_G(B_2)} \in \{4,6\}$. However, if $\abs{N_G(B_2)} = 4$, then $\abs{G} = \abs{N_G(B_3)}$, making the configuration axial, a contradiction. So, $\abs{N_G(B_2)} = 6$, and then we have  $\abs{N_G(B_3)} \in \{6,8,10\}$. In summary,
\[\abs{N_G(B_1)} = 4,\abs{N_G(B_2)} = 6,\abs{N_G(B_3)} \in \{6,8,10\}.\]

\begin{claim}
If $(r,s) = (0,3)$ then $\abs{N_G(B_3)}\neq6$.
\end{claim}
\begin{proofclaim}
If $\abs{N_G(B_3)} = 6$, then $\abs{G} = 12$. In particular,  $C_G(B_3) \triangleleft N_G(B_3) \triangleleft G$, implying $C_G(B_3)$ is a normal Sylow $3$-subgroup of $G$. But, since $\abs{N_G(B_2)} = 6$ as well, the same is true of $C_G(B_2)$, a contradiction.
\end{proofclaim}

\begin{claim}\label{cl.rs03n8}
If $(r,s) = (0,3)$ and $\abs{N_G(B_3)}=8$, then $G\cong \Sym(4)$.
\end{claim}
\begin{proofclaim}
In this case, $\abs{G} = 24$, and the orbit $\mathcal{O}$ containing $B_2$ has size $4$. The kernel $K$ of the action on $\mathcal{O}$ has size dividing $\abs{N_G(B_2)} = 6$ and also dividing $4$ by Remark~\ref{rem.BasicPrinciples}.
We rule out $\abs{K}=2$.
If $\abs{K} = 2$ and $k\in K^*$, then $k$ is central in $G$, so $G$ normalizes the axis of $k$, contradicting that $G$ is nonaxial. We conclude that $G$ acts faithfully on $\mathcal{O}$, so $G\cong \Sym(4)$.
\end{proofclaim}

\begin{claim}\label{cl.rs03n10}
If $(r,s) = (0,3)$ and $\abs{N_G(B_3)}=10$, then $G\cong \Alt(5)$.
\end{claim}
\begin{proofclaim}
Here $\abs{G} = 60$. Also notice that centralizers of axes have order $2$, $3$, or $5$, implying  every nontrivial element of $G$ has order $2$, $3$, or $5$. 

We first claim that distinct Sylow $2$-subgroups intersect trivially. So, suppose the intersection of  Sylow $2$-subgroups $P$ and $Q$ contains some nontrivial $g$. As $P$ is abelian, its elements normalize the axis of $g$ and similarly for $Q$; thus $P,Q \le N_G(A_g)$.  Since $g$ has order $2$, $A_g$ is in the orbit of $B_1$, so $N_G(A_g)$ has order $4$, implying  $P = N_G(A_g) =  Q$ as desired.
    
We now find the number of Sylow $2$-subgroups by counting the set $I$ of involutions.
Every involution must have an axis in the orbit $\mathcal{O}$ of $B_1$, so
\[ \abs{I} = \sum_{A \in \mathcal{O}} \abs{C^*_G(A)} = \abs{C^*_G(B_1)} \cdot \abs{G:N_G(B_1)} = 1 \cdot 15.\]
As Sylow $2$-subgroups intersect trivially, $G$ must has five. Letting $G$ act on its Sylow $2$-subgroups, routine analysis---expedited by the additional information we have---shows that the action is faithful, so $G\cong \Alt(5)$.
%
\end{proofclaim}

This completes the classification of the nonaxial 3D-rotation groups. 


\begin{claim}
For each $A\in \mathbb{A}$, $C_G(A)$ is cyclic. 
\end{claim}
\begin{proofclaim}
This is mostly immediate from our classification as $C_G(A)$  has prime order except when $G=\Sym(4)$ and $A$ is conjugate to $B_3$. In this exceptional case, only centralizers of axes conjugate to $B_3$ are large enough to accommodate a $4$-cycle, and as such a centralizer has size exactly $4$, it is cyclic, as desired.
\end{proofclaim}
\end{proof}\setcounter{claim}{0} 

\section{Recognition of icosahedral modules}\label{S.Alt5Recognition}

This section is devoted to the proof of Theorem~\ref{thm.Alt5Recognition}.
Our methods are  elementary and require little background. In \S\S~\ref{s.trace-adjoint}, we discuss basic properties of the \emph{trace} and \emph{adjoint} maps associated to endomorphisms of an abelian group; these feature frequently in our proof. The brief \S\S~\ref{s.setup} establishes notation and  outlines our strategy. The key ingredients of our proof are developed in \S\S~\ref{s.class-sums}--\ref{s.local} and then assembled in \S\S~\ref{s.recognition-proof}.

\subsection{Trace and adjoint maps}\label{s.trace-adjoint}

\begin{definition*}
Let $G$ be a group and $V$ a $G$-module. If $g\in G$ has finite order $m$, we define maps $\ad_g,\tr_g\in \End(V)$ as follows:
\begin{itemize}
    \item $\ad_g = 1-g$ is called the \textbf{adjoint} map;
    \item $\tr_g = 1+g+\cdots+g^{m-1}$ is called the \textbf{trace} map.
\end{itemize}
When the module $V$ is clear from the context, we use notation $B_g = \ad_g(V)$ and $C_g = \tr_g(V)$ for the images of these maps. 
\end{definition*}


\begin{fact}[{\cite[Lemma~2.1]{CDWY24}}]\label{fact.AdjointTrace}
Let $G$ be a group and $V$ a $G$-module. If $g\in G$ has finite order, then the following hold.
\begin{enumerate}
    \item\label{fact.AdjointTrace.i:CCentral} $\ad_g\circ \tr_g = 0$, so $C_g$ is centralized by $g$.
    \item\label{fact.AdjointTrace.i:traceonB} $\tr_g\circ \ad_g = 0$, so if $|g| = 2$, then $B_g$ is inverted by $g$.
    \item\label{fact.AdjointTrace.i:BCIntersection} Every element of $B_g \cap C_V(g)$ has order dividing $|g|$.
\end{enumerate}
\end{fact}

\begin{lemma}[{cf.~\cite[Coprimality Lemma]{CDW23}}]\label{lem.m-divisible-decomposition}
Let $G$ be a group, $V$ a $G$-module, and $g\in G$ of finite order $m$. If $V$ is $m$-divisible, then $V = B_g + C_g$, and if $V$ is uniquely $m$-divisible, then the sum is direct. 
\end{lemma}
\begin{proof}
    Assume $V$ is $m$-divisible. Let $v\in V$; write $v = mw$ for some $w\in V$. Set $c = \tr_g(w)$ and $b = \sum_{i=1}^{m-1} (w - g^i w)$. Notice that $b + c = mw = v$. Also, $b\in  B_g$ since $w - g^i w =  \ad_g\left(\sum_{j = 1}^i g^{j-1} w \right)$. Thus $V = B_g + C_g$.
And if $V$ is uniquely $m$-divisible, then $V$ has no nontrivial elements of order dividing $m$, so Fact~\ref{fact.AdjointTrace}(\ref{fact.AdjointTrace.i:CCentral},\ref{fact.AdjointTrace.i:BCIntersection}) implies $B_g \cap C_g = 0$.
\end{proof}

\begin{remark}\label{rem.C-equals-Centralizer}
One consequence of Lemma~\ref{lem.m-divisible-decomposition} is that unique $m$-divisibility implies $C_g = C_V(g)$ when $|g| = m$.  
\end{remark}

\begin{lemma}[{cf.~\cite[Weight Lemma]{CDW23}}]\label{lem.2-divisible-Klein-decomposition}
Let $K=\{1, \alpha_1, \alpha_2, \alpha_3\}$ be a Klein $4$-group and  $V$ a $K$-module. Define the following subgroups of $V$.
\begin{align*}
V_{\alpha_1} & = \tr_{\alpha_1} \circ \ad_{\alpha_2} \circ \ad_{\alpha_3}(V) & V_{\alpha_3} &= \ad_{\alpha_1} \circ \ad_{\alpha_2} \circ \tr_{\alpha_3}(V)\\
V_{\alpha_2} &= \ad_{\alpha_1} \circ \tr_{\alpha_2} \circ \ad_{\alpha_3}(V) & C_K &= \tr_{\alpha_1} \circ \tr_{\alpha_2} \circ \tr_{\alpha_3}(V)
\end{align*}
If $V$ is $2$-divisible, then $V = V_{\alpha_1} + V_{\alpha_2} + V_{\alpha_3} + C_K$, and if $V$ is uniquely $2$-divisible, then the sum is direct.
\end{lemma}
\begin{proof}
Follow the argument of \cite[Lemma~2.14]{CDWY24} using Lemma~\ref{lem.m-divisible-decomposition} in place of their Fact~2.12.
\end{proof}

\subsection{Setup and strategy}\label{s.setup}

We adopt the following setup for the remainder of the present section. 

\begin{setup}
Let $G = \Alt(5)$ and $V$ be a uniquely $2$-divisible $G$-module  that either has prime characteristic or is divisible. In the case that $V$ is $p$-divisible and $g\in \Alt(5)$ has order $p$, we assume  $\tr_g(V)$ is inverted by each element  of $N_{G}(\langle g \rangle) \setminus \langle g \rangle$; we will refer to this property (rather vaguely) as the \emph{normalizer assumption}.
\begin{itemize}
    \item For $g\in G$, we write $C_g = \tr_g(V)$ and $B_g = \ad_g(V)$. 
    \item Let $\mathcal{C}_2,\mathcal{C}_3,\mathcal{C}_{5,1}, \mathcal{C}_{5,2}$ denote the conjugacy classes of elements of order $2$, $3$, and $5$.
    Then let $c_2,c_3,c_{5,1},c_{5,2} \in \Z[G]$ denote the corresponding class sums (e.g.~$c_2 =\sum_{g\in \mathcal{C}_2} g$). We also set $c_5 = c_{5,1} + c_{5,2}$.
\end{itemize}
\end{setup}

\begin{remark}\label{rem.NormalizerAssumption}
The structure of $\Alt(5)$ implies that the normalizer assumption is equivalent to the following statement: for all $\alpha \in N_{G}(\langle g \rangle)\setminus \langle g \rangle$, $\tr_\alpha \circ \tr_g = 0$. We produce additional characterizations in Lemma~\ref{lem.A5.dihedral-sum}.
\end{remark}

Let us outline our approach to showing $V$ is icosahedral (as defined in \S\S~\ref{s.results}) with respect to some generating set $\{\alpha,\gamma\}$ for $G$ with $\abs{\alpha} = 2$ and $\abs{\gamma} = 3$. Identifying the subgroups $Re_1,Re_2,Re_3$ is straightforward as the definition of an icosahedral module forces  $Re_3 = C_\alpha$, 
then implying $Re_2 = \gamma C_\alpha$ and $Re_3 = \gamma^{2} C_\alpha$. Next we need a candidate for $\phi$.  Proposition~\ref{prop.A5.class-sums} will show that $\frac{1}{4}c_{5,1}$ and $\frac{1}{4}c_{5,2}$ have the desired property, and our choice of the two depends only on our generating set for $\Alt(5)$. For example, if $\alpha\gamma \in \mathcal{C}_{5,1}$, we take $\phi = \frac{1}{4}c_{5,1}$.
It then remains to prove the remaining relations for $\alpha$. The condition $\alpha e_2 = -e_2 + \phi e_3$ can be recast as $(\alpha\gamma + \gamma\alpha) e_3 = \phi e_3$, and we will  show in Proposition~\ref{prop.A5.local} that  $(\alpha\gamma + \gamma\alpha)$ and $\frac{1}{4}c_{5,1}$ agree on  $C_\alpha$. The relation  $\alpha e_1 = -e_1 + \phi e_3$ is similarly addressed, this time by showing $(\alpha\gamma^{-1} + \gamma^{-1}\alpha)$ and $\frac{1}{4}c_{5,1}$ agree on  $C_\alpha$.

\subsection{The class sums}\label{s.class-sums}

The present goal is to prove the following.

\begin{proposition}\label{prop.A5.class-sums}
In $\End(V)$, we have that 
\[c_2 = -5, c_3 = 0, c_{5} = 4,\] and both $c_{5,1}$ and $c_{5,2}$ satisfy $x^2 = 4(x+4)$.
\end{proposition}

We  record a few  basic observations.

\begin{lemma}\label{lem.Alt(5)-action-Klein}
Let $K=\{1, \alpha_1, \alpha_2, \alpha_3\}$ be a Klein $4$-subgroup of $G$. Then $V = C_{\alpha_1} + C_{\alpha_2} + C_{\alpha_3}$ and $C_V(K) = 0$. In particular, $C_V(G) = 0$.
\end{lemma}
\begin{proof}
We adapt the argument of \cite[Lemma~4.14]{CDWY24} to show $C_V(K) = 0$.

In the notation of Lemma~\ref{lem.2-divisible-Klein-decomposition}, we have $V = V_{\alpha_1} \oplus V_{\alpha_2} \oplus V_{\alpha_3} \oplus C_K$. Since $C_{K}\le C_{\alpha_1}\cap C_{\alpha_2}$, $C_{K}$ is centralized by $\alpha_1$ and, by the normalizer assumption,  also inverted by $\alpha_1$. So as $V$ has no $2$-torsion, $C_K=0$. In particular, since $V_{\alpha_i} \le C_{\alpha_i}$, we certainly have $V = C_{\alpha_1} + C_{\alpha_2} + C_{\alpha_3}$.

Let $w\in C_V(K)$. We  now run the computation in \cite[Lemma~4.14]{CDWY24} to see that $w = 0$. Write $w = v_1 + v_2 + v_3$ with $v_i \in V_{\alpha_i}$. Applying $\alpha_1$ (and remembering that $w\in C_V(K)$), we get $v_1 + v_2 + v_3 = v_1 - v_2 - v_3$, implying $2(v_2 + v_3) = 0$ so $v_2 + v_3 = 0$. Similarly, $v_1 + v_3 = 0 = v_1 + v_2$, and combining, we find that $w = 0$.
\end{proof}

\begin{lemma}\label{lem.sum-G-in-End(V)}
In $\End(V)$ we have that $\sum_{g\in G} g = 0.$
\end{lemma}
\begin{proof}
Let $z = \sum_{g\in G} g$; then for each $g\in G$, $gz = z$. Thus, for all $v\in V$, $g(zv) = zv$, so $zv \in C_V(G)$, implying, by Lemma~\ref{lem.Alt(5)-action-Klein}, that $zv = 0$.
\end{proof}

\begin{lemma}\label{lem.A5.dihedral-sum}
Let $g\in G$ be nontrivial, and write $D = N_{G}(\langle g \rangle)$. If $|g|=p$ and $V$ is $p$-divisible, then for all $\alpha \in D\setminus \langle g \rangle$, $\tr_\alpha \circ \tr_g = \tr_g \circ \tr_\alpha = \sum_{d\in D} d=0$.
\end{lemma}
\begin{proof}
The  point is that $D$ is dihedral. For $\alpha \in D\setminus \langle g \rangle$, $\tr_\alpha = 1+\alpha$, and the normalizer assumption yields $\tr_\alpha \circ \tr_g = 0$. Moreover, 
\[\tr_\alpha \tr_g  = \sum_{g\in \langle g\rangle} (g + \alpha g) = \sum_{d\in D} d = \sum_{g\in \langle g\rangle} (g + g\alpha )   =  \tr_g\tr_\alpha.\qedhere\]
\end{proof}

\begin{proof}[Proof of Proposition~\ref{prop.A5.class-sums}]
Lemma~\ref{lem.A5.dihedral-sum} is our main tool to compute the $c_i$. The lemma requires $p$-divisibility, which we only know for $p= 2$. However, our assumptions  imply that $V$ will be $p$-divisible for all except possibly one prime $p$. This  allows us to  directly compute $c_2$ and one of $c_3, c_5$; the remaining one is then  found using Lemma~\ref{lem.sum-G-in-End(V)}.

We first compute $c_2$. Let $D_1,\ldots, D_5$ be the Sylow $2$-subgroups of $G$, and note that $\mathcal{C}_2 =\bigcup_i D_i^*$. Each $D_i$ is the normalizer of an order $2$ subgroup, so applying Lemma~\ref{lem.A5.dihedral-sum}, we find that
\[ c_2 = \sum_{i=1}^5 \sum_{d\in D_i^*} d = \sum_{i=1}^5\left( -1+\sum_{d\in D_i} d\right) = \sum_{i=1}^5\left(  -1 +0\right) = -5.\]

Now assume $V$ is $3$-divisible, and let $D_1,\ldots, D_{10}$ denote the normalizers of the Sylow $3$-subgroups. Each $3$-cycle lies in a unique $D_i$, and each involution of $G$ appears in precisely two:  $(ab)(cd)$ normalizes the Sylow $3$-subgroups $\langle (abe) \rangle$ and $\langle (cde) \rangle$ and no others. This shows that $\bigcup_i D_i^*$, viewed as a multiset, is equal to $\mathcal{C}_3\cup 2\mathcal{C}_2$, so again using Lemma~\ref{lem.A5.dihedral-sum}, we find that
\[ c_3 + 2c_2 = \sum_{i=1}^{10} \sum_{d\in D_i^*} d = \sum_{i=1}^{10}\left( -1+\sum_{d\in D_i} d\right) = -10.\]

And if $V$ is $5$-divisible, we work with the normalizers $D_i$ of the Sylow $5$-subgroups. Noting that $5$-cycles are in a unique $D_i$ and involutions are in exactly two, we see that $\bigcup_i D_i^* = \mathcal{C}_5\cup 2\mathcal{C}_2$ and 
\[ c_5 + 2c_2 = \sum_{i=1}^{6} \sum_{d\in D_i^*} d = \sum_{i=1}^{6}\left( -1+\sum_{d\in D_i} d\right) = -6.\]

In total, we know $c_2 = -5$ and also the value of one of $c_3$ or $c_5$. The remaining value is found using $1+c_2 + c_3 + c_5 = 0$ from Lemma~\ref{lem.sum-G-in-End(V)}. 

The condition on $c_{5,i}$ now follows from the fact that $c_{5,1}^2 = 12 + 3c_3 +  c_{5}+ 4c_{5,1}$ and similarly for $c_{5,2}$. This holds in $\Z[G]$ so can be seen various ways such as by direct computation (by hand or with GAP \cite{GAP4}) or via the character-theoretic formula for multiplying in $Z(\mathbb{C}[G])$ with the class sums as a basis (see for example \cite[Problem~3.9]{IsM94}). 
\end{proof}

\subsection{Local Equations}\label{s.local}

\begin{notation}
Fix a generating set $\{\alpha,\gamma\}$ for  $G=\Alt(5)$ with $|\alpha| = 2$ and $|\gamma| = 3$. 
Fix additional notation as follows:
\begin{itemize}
\item set $\sigma_1 = \alpha\gamma$ and $\sigma_2 = \gamma\alpha$ (which are $\alpha$-conjugate $5$-cycles);
\item label the conjugacy classes so that $\sigma_1,\sigma_2 \in \mathcal{C}_{5,1}$;
\item if $e_1,e_2\in \End(V)$, we write $e_1\equiv e_2$ to mean $e_1\tr_\alpha = e_2 \tr_\alpha$. 
\end{itemize}
For concreteness, one may take $\alpha = (12)(34)$ and $\gamma = (135)$ in which case $\sigma_1 = (14352)$ and $\sigma_2 = (12345)$.
\end{notation}

\begin{remark*}
Since $V$ is $2$-divisible and $C_\alpha \le C_V(\alpha)$, $e_1\equiv e_2$ may alternatively be read as $e_1$ and $e_2$ agree on $C_\alpha$.
\end{remark*}

We aim to prove the following local equation.

\begin{proposition}\label{prop.A5.local}
$\sigma_1 + \sigma_2 \equiv  \frac{1}{4}c_{5,1}  \equiv \sigma_1^{-1} + \sigma_2^{-1}$.
\end{proposition}
\begin{proof}
We  identify the $\alpha$-orbits on $\mathcal{C}_{5,1}$ in terms of $\sigma_1$ and $\sigma_2$ and  study how the sums over these orbits---which includes $\sigma_1 + \sigma_2$---act on  $C_\alpha$. 

For $\beta\in C_G(\alpha)$, $\beta$-conjugates of sums over $\alpha$-orbits will be equivalent on $C_\alpha$ (see Claim~\ref{claim-prop.A5.local-basic-equalities} below), so we also identify the  $\beta$-orbits on  $\mathcal{C}_{5,1}$ for $\beta$ chosen as follows: let $\beta'$  be the unique involution of $C_G(\alpha)$ that inverts $\gamma$ and then set $\beta = \alpha\beta'$. Thus, $\beta$ inverts both $\sigma_1$ and $\sigma_2$.
 The elements of $\mathcal{C}_{5,1}$ are shown in Figure~\ref{fig:sigma1-class} with conjugacy via $\alpha$ and $\beta$ denoted via arrows. 

\begin{figure}[ht]
    \centering
    \begin{tikzpicture}[yscale = 0.75]
        \footnotesize
        \node (Q1) at (-5,1) {$\sigma_1$};
        \node (Q2) at (-3,1) {$\sigma_2$};
        \node (Q3) at (-5,-1) {$\sigma_1^{-1}$};
        \node (Q4) at (-3, -1) {$\sigma_2^{-1}$};

        \node (Q5) at (-1,1) {$\sigma_2^2\alpha$};
        \node (Q6) at (1,1) {$\sigma_1^2\alpha$};
        \node (Q7) at (-1,-1) {$\sigma_2^3\alpha$};
        \node (Q8) at (1,-1) {$\sigma_1^3\alpha$};

        \node (Q9) at (3,1) {$\sigma_2^{-1}\sigma_1$};
        \node (Q10) at (5,1) {$\sigma_1^{-1}\sigma_2$};
        \node (Q11) at (3,-1) {$\sigma_2\sigma_1^{-1}$};
        \node (Q12) at (5,-1) {$\sigma_1\sigma_2^{-1}$};

        \foreach \i / \j in {1/2,5/6,9/10}  {
            \draw[red, thick, <->] (Q\i) -- (Q\j) node[midway, below] {$\alpha$};}
          \foreach \i / \j in {3/4,7/8,11/12}  {
            \draw[red, thick, <->] (Q\i) -- (Q\j) node[midway, above] {$\alpha$};}

        \foreach \i / \j in {1/3,5/7,9/11}  {
            \draw[blue, thick, <->] (Q\i) -- (Q\j) node[midway, right] {$\beta$};}
          \foreach \i / \j in {2/4,6/8,10/12}  {
            \draw[blue, thick, <->] (Q\i) -- (Q\j) node[midway, left] {$\beta$};}

    \end{tikzpicture}
    \caption{Conjugacy class of $\sigma_1$: $\mathcal{C}_{5,1}$}
    \label{fig:sigma1-class}
\end{figure}

\begin{notation}
We name three sums of $\alpha$-conjugates;  Claim~\ref{claim-prop.A5.local-basic-equalities} will show that every sum of $\alpha$-conjugates in $\mathcal{C}_{5,1}$ is equivalent to one of these. 
\begin{align*}
\phi &= \sigma_1 + \sigma_2 & u &= \sigma_1^2\alpha + \sigma_2^2\alpha & r &=  \sigma_2^{-1}\sigma_1 +  \sigma_1^{-1}\sigma_2  
\end{align*}
As we work towards proving that $4\phi \equiv  c_{5,1}$, it may help to keep the classical setting in mind: there, $\phi$ acts on $C_\alpha$ as the golden ratio or its conjugate, $u\equiv 1$, and  $\phi - r\equiv 1$ (meaning $r$ acts as the negative of the conjugate of $\phi$).
\end{notation}

\begin{claim}\label{claim-prop.A5.local-basic-equalities}
For all $e \in \End(V)$, one has that $e \equiv e\alpha$ and $e + \alpha e \alpha \equiv \beta (e + \alpha e \alpha) \beta$. In particular, \[c_{5, 1} \equiv 2(\phi + u + r).\]
\end{claim}
\begin{proofclaim}
Since $\ad_\alpha\tr_\alpha = 0$, we have $e(1-\alpha)\tr_\alpha = 0$, so  $e \equiv e\alpha$. 

Now, let $f = \alpha e \alpha$, and notice that $e + f$ commutes with $\tr_\alpha$. Since  $V$ is $2$-divisible, the normalizer assumption applied to $\langle \alpha \rangle$ tells us that  $\beta \tr_\alpha = -\tr_\alpha$ and  also implies that $\tr_\alpha \tr_\beta = \tr_\beta \tr_\alpha = 0$ (see Lemma~\ref{lem.A5.dihedral-sum}). Therefore,
\[\beta(e + f)\beta \tr_\alpha = - \beta (e+f)\tr_\alpha = - \beta \tr_\alpha (e+f) = \tr_\alpha (e+f) = (e+f) \tr_\alpha,\]
so $\beta (e+f) \beta \equiv e+f$, as desired.

The local equation for $c_{5,1}$ now follows since $\phi+u+r$ may be written in the form  $e + \alpha e \alpha$ and $c_{5,1} = (\phi + u + r) + \beta(\phi + u + r)\beta$.
\end{proofclaim}

\begin{claim}\label{claim-prop.A5.phisquared}
$\phi^2 \equiv u + \phi \equiv 2 + r$. Moreover, $u\phi \equiv \phi u \equiv u + r$ (which we will only use in characteristic $3$).
\end{claim}
\begin{proofclaim}
We use Claim~\ref{claim-prop.A5.local-basic-equalities} repeatedly.

We have $\phi^2 = (\sigma_1 + \sigma_2)^2 = \sigma_1^2 + \sigma_2^2 + \alpha\gamma^2\alpha + \gamma^2 \equiv u + \alpha\gamma^2 + \gamma^2\alpha = u + \sigma_2^{-1} + \sigma_1^{-1}
\equiv u + \phi$. Also, $\phi^2 \equiv (\sigma_1 + \sigma_2)(\sigma_1^{-1} + \sigma_2^{-1})= 2 + \sigma_2\sigma_1^{-1} + \sigma_1\sigma_2^{-1} \equiv 2 + r,$ establishing the main statement of the claim.

Next, observe that
$u\phi = (\sigma_1^2 + \sigma_2^2)\alpha (\sigma_1 + \sigma_2) =  (\sigma_1^2 + \sigma_2^2)(\sigma_1 + \sigma_2) \alpha = (\sigma_1^3 + \sigma_2^3 +\sigma_1^2\sigma_2 + \sigma_2^2\sigma_1 )\alpha \equiv u + \sigma_1^2\sigma_2 + \sigma_2^2\sigma_1$. Additionally,
$\sigma_1^2\sigma_2 + \sigma_2^2 \sigma_1 =
(\alpha\gamma)^2(\gamma\alpha) + (\gamma\alpha)^2(\alpha\gamma)
 \equiv \alpha \gamma\alpha \gamma^2 + \gamma\alpha \gamma^2\alpha
= \sigma_1\sigma_2^{-1} + \sigma_2 \sigma_1^{-1} \equiv r$. In total, $u\phi \equiv u + r$, and  
analogous computations show $\phi u \equiv u + r$ as well. 
\end{proofclaim}

We are now ready to show that $u \equiv 1$. Here our methods differ dramatically depending on whether $V$ is $3$-divisible or not.

\begin{notation}
We set $\delta = \sigma_1^{-2}\sigma_2^2$. Note that $\delta$ is a $3$-cycle inverted by $\alpha$:  taking $\alpha = (12)(34)$ and $\gamma = (135)$, we have $\delta = (345)$. Moreover, $\beta' = (13)(24)$ (as defined earlier), so  $\beta' = \sigma_1^{2}\delta$.
\end{notation}

\begin{claim}\label{claim-prop.A5.3-divisible}
If $V$ is $3$-divisible, then $u \equiv 1$.
\end{claim}
\begin{proofclaim}
By Lemma~\ref{lem.A5.dihedral-sum}, $\tr_\delta \tr_\alpha = 0$ and $\tr_{\beta'} \tr_\alpha = 0$, so $1+\delta \equiv -\delta^2$ and $\beta'\equiv -1$. Thus, 
\[u \equiv \sigma_1^2 + \sigma_2^2 = \sigma_1^2 (1 + \delta) \equiv - \sigma_1^2 \delta^2 = - \beta' \equiv 1.\qedhere\]
\end{proofclaim}

It remains to consider when $V$ is not $3$-divisible, which (by our assumptions) implies that $\charac V = 3$.

\begin{claim}\label{claim-prop.A5.char3}
If $\charac V = 3$, then $u \equiv 1$.
\end{claim}
\begin{proofclaim}
In characteristic $3$, Claims~\ref{claim-prop.A5.local-basic-equalities} and \ref{claim-prop.A5.phisquared} imply $c_{5, 1} \equiv \phi + u - 1$. Also, by Proposition~\ref{prop.A5.class-sums}, $c_{5, 1}$ satisfies $x^2 = 4(x + 4) = x + 1$, so again using Claim~\ref{claim-prop.A5.phisquared}, we find that
\begin{align*}
0 & = (\phi + u - 1)^2 - (\phi + u - 1) - 1\\
& = \phi^2 + u^2 + u\phi + \phi u + 1\\
& \equiv \phi + u + u^2 +2(2 u + \phi - 2) + 1\\
& = u^2 - u.
\end{align*}
Thus, $u(u-1) \equiv 0$, so to prove the claim, it suffices to show that the restriction of $u$  to $C_\alpha$ is injective. 

Let $v\in C_\alpha \cap \ker u$. As in the proof of Claim~\ref{claim-prop.A5.3-divisible}, $u \equiv \sigma_1^2 + \sigma_2^2 = \sigma_1^2 (1 + \delta)$, so $(1 + \delta) v = \sigma_1^{-2} u v = 0$. Thus $\delta v = -v$.  But $\delta$ is a $3$-cycle, so $v = \delta^3 v = - v$. As $V$ has no $2$-torsion, $v=0$, and $u$ is injective.
\end{proofclaim}

\begin{claim}
$4\phi \equiv  c_{5,1}$.
\end{claim}
\begin{proofclaim}
With $u \equiv 1$ in hand, we first revisit Claim~\ref{claim-prop.A5.phisquared} to get $r \equiv \phi -1$, and then return to Claim~\ref{claim-prop.A5.local-basic-equalities}:
\[c_{5,1} \equiv 2\phi + 2u + 2r \equiv 2\phi + 2 + 2(\phi - 1) \equiv 4 \phi.\qedhere\]
\end{proofclaim}
As $\phi = \sigma_1 + \sigma_2 \equiv \sigma_1^{-1} + \sigma_2^{-1}$,  we are done.
\end{proof}\setcounter{claim}{0} 

\subsection{Proof of Theorem~\ref{thm.Alt5Recognition}}\label{s.recognition-proof}

We first recall the definition of an icosahedral module and then proceed to identitfy our $\Alt(5)$-module $V$.

\begin{definition*}[Icosahedral Modules]
Let $R = \Z[\phi]$ with $\phi$ satisfying $\phi^2 = \phi +1$. Fix a generating set $\{\alpha,\gamma\}$ for  $\Alt(5)$ with $|\alpha| = 2$,  $|\gamma| = 3$, and $|\alpha\gamma| = 5$. Define $\icosa(R)$ (with respect to $\alpha,\gamma$) to be the $R[\Alt(5)]$-module $Re_1 \oplus Re_2 \oplus Re_3$ with action given as follows:
\begin{align*}
\alpha e_1 & = -e_1 + \phi e_3 & \gamma e_1 &= e_3\\
\alpha e_2 & = -e_2 + \phi e_3 & \gamma e_2 &= e_1\\
\alpha e_3 & = e_3 & \gamma e_3 &= e_2
\end{align*}
Then, for any  $R$-module $L$ without $2$-torsion, set $\icosa(L) = \icosa(R) \otimes_R L$ with $\Alt(5)$ acting trivially on $L$.
\end{definition*}

\begin{proof}[Proof of Theorem~\ref{thm.Alt5Recognition}]
Choose an involution  $\alpha\in \Alt(5)$ and a conjugacy class of $5$-cycles $\mathcal{C}$. 
From this, choose a $3$-cycle $\gamma$ such that $\alpha\gamma \in \mathcal{C}$ and note that $\Alt(5) = \langle \alpha,\gamma\rangle$. 
Adopt the initial notation of \S\S~\ref{s.local}; in particular: $ \sigma_1 = \alpha\gamma$, $ \sigma_2 = \gamma\alpha$, and $\mathcal{C}_{5,1} = \mathcal{C}$. Define $\phi$ to be the image of $\frac{1}{4}c_{5,1}$ in $\End(V)$; by Proposition~\ref{prop.A5.class-sums}, $\phi^2 = \phi + 1$. Additionally, since $\phi$ commutes with $\Alt(5)$, $\phi$ acts on each $C_g$ for $g\in \Alt(5)$ so may  be viewed as an element of $\End(C_g)$.

Define $R = \Z[\phi]$. Set $L= C_\alpha$ considered as an $R$-module---with action of $\phi$ defined as above---and also as a trivial $\Alt(5)$-module. Construct $\icosa(L)$  using this $L$ and the above generating set; we will show that $\icosa(L) \cong V$ as $\Alt(5)$-modules. Note that while $\Alt(5)$ acts trivially on $L$ in the construction of  $\icosa(L)$, $\Alt(5)$ most certainly does not act trivially on $C_\alpha$ as a subgroup of $V$.

Consider the $R$-module morphism $f:\icosa(L) \rightarrow V$ defined via \[f(e_i\otimes \ell) = \gamma^{-i}\ell.\]

\begin{claim}\label{claim-thm.Alt5Recognition.morphism}
$f$ is an $\Alt(5)$-module morphism.
\end{claim}
\begin{proofclaim}
By design, $f$ commutes with the action of $\gamma$, and by choice of $L$, $f(\alpha(e_3\otimes \ell)) = \ell = \alpha f(e_3\otimes \ell)$. 
Now, using Proposition~\ref{prop.A5.local} and again using that $\alpha$ centralizes $L$, we have
\begin{align*}
\alpha f(e_2\otimes \ell) 
 = \alpha\gamma \ell 
& = \alpha\gamma \ell + \gamma \alpha\ell - \gamma \alpha\ell \\
& = (\sigma_1 + \sigma_2)\ell -\gamma\ell \\
& = \phi\ell -\gamma \ell \\
& = f((\phi e_3 - e_2)\otimes \ell)\\
& = f(\alpha(e_2\otimes \ell)).
\end{align*}
Similarly, $\alpha f(e_1\otimes \ell) = (\sigma^{-1}_1 + \sigma^{-1}_2)\ell -\gamma^{-1}\ell = \phi\ell -\gamma^{-1} \ell = f(\alpha(e_1\otimes \ell)).$
\end{proofclaim}

\begin{claim}
$f$ is surjective.
\end{claim}
\begin{proofclaim}
Let $W$ be the image of $f$. It certainly contains $C_\alpha$ and is $\Alt(5)$-invariant by Claim~\ref{claim-thm.Alt5Recognition.morphism}. Thus, since involutions are conjugate $\Alt(5)$, $C_\beta \le W$ for every involution $\beta$, implying that $W = V$ by Lemma~\ref{lem.Alt(5)-action-Klein}.
\end{proofclaim}

\begin{claim}
$f$ is injective.
\end{claim}
\begin{proofclaim}
Suppose $w = e_1\otimes \ell_1 + e_2\otimes \ell_2 + e_3\otimes \ell_3$ is in the kernel. By Claim~\ref{claim-thm.Alt5Recognition.morphism},  $\tr_\alpha(w), \tr_\alpha (\gamma^{-1} w), \tr_\alpha (\gamma w)$ are also in the kernel, which produces the following system in $C_\alpha$:
\begin{align*}
\phi\ell_1  +  \phi\ell_2  +  2\ell_3 &= 0\\
\phi\ell_1  +  2\ell_2  +  \phi\ell_3 &= 0\\
2\ell_1  +  \phi\ell_2  +  \phi\ell_3 &= 0.
\end{align*}
Using that both $2$ and $\phi$ are invertible in $\End(V)$ (the first by unique $2$-divisibility and the second as $\phi^2 = \phi + 1$ implies $\phi(\phi - 1)=1$), the system easily yields $\ell_1 = \ell_2 = \ell_3 = 0$.
\end{proofclaim}
\end{proof}\setcounter{claim}{0}

\section{The corollaries}\label{S.Corollaries}

We now address the corollaries, which we restate here for convenience.

\begin{namedthm}{Corollary~\ref{cor.Alt5Recognition-3Drot}}
Let $V$ be a nonaxial 3D-rotation module for $\Alt(5)$ with respect to the axes $A_g = \tr_g(V)$. If $V$ is uniquely $2$-divisible and  either has prime characteristic or is divisible, then the conclusion of Theorem~\ref{thm.Alt5Recognition} holds.
\end{namedthm}

\begin{namedthm}{Corollary~\ref{cor.Alt5Recognition-dim}}
Let $V$ be a faithful $\Alt(5)$-module with an additive dimension. Assume $V$ is dim-connected and $\dim V = 3$. If $V$ is without $2$-torsion, then  the conclusion of Theorem~\ref{thm.Alt5Recognition} holds and $\dim L = 1$.
\end{namedthm}

\subsection{Proof of Corollary~\ref{cor.Alt5Recognition-3Drot}}

Corollary~\ref{cor.Alt5Recognition-3Drot} follows almost immediately from Theorems~\ref{thm.ClassificationRotationGroups} and \ref{thm.Alt5Recognition}.

\begin{proof}[Proof of Corollary~\ref{cor.Alt5Recognition-3Drot}]
Take $V$ to be a nonaxial 3D-rotation module for $G=\Alt(5)$ with respect to the choice of axes $A_g = \tr_g(V)$. Also assume that $V$ is uniquely $2$-divisible and  either has prime characteristic or is divisible. To apply Theorem~\ref{thm.Alt5Recognition}, it only remains to verify the normalizer assumption. 

By Theorem~\ref{thm.ClassificationRotationGroups}, $C_G(A_g)$ is cyclic, so by the structure of $\Alt(5)$, $C_G(A_g) = \langle g\rangle$. The definition of a  3D-rotation module implies $N_G(A_g) = N_G(C_G(A_g))$, so $N_G(A_g) = N_G(\langle g\rangle)$. As such, the normalizer assumption of Theorem~\ref{thm.Alt5Recognition} follows from the corresponding assumption for a 3D-rotation module.
\end{proof}

\subsection{Proof of Corollary~\ref{cor.Alt5Recognition-dim}}\label{s.Alt5Recognition-dim}

Corollary~\ref{cor.Alt5Recognition-dim} is an application of Theorem~\ref{thm.Alt5Recognition} to modules that carry a well-behaved notion of dimension (which might not be associated to a vector space structure). We work in the context of \emph{dim-connected modules with an additive dimension} as presented in~\cite{CDW23}. Here we collect a handful of basic definitions and results for these modules; for more details, we refer to~\cite{CDW23} and \cite{CDWY24}.

\begin{definition*}
A \textbf{modular universe} $\cU$ is a collection of abelian groups $\Ob(\cU)$ and homomorphisms $\Ar(\cU)$ between them that satisfy the following closure properties.
\begin{itemize}
\item{}\textsc{[inverses]}
If $f \in \Ar(\cU)$ is an isomorphism, then $f^{-1} \in \Ar(\cU)$.
\item{}\textsc{[products]} If $V_1, V_2 \in \Ob(\cU)$ and $f_1, f_2 \in \Ar(\cU)$, then 
$V_1 \times V_2 \in \Ob(\cU)$, and 
	$\Ar(\cU)$ contains $f_1 \times f_2$,  the projections $\pi_i : V_1\times V_2 \rightarrow V_i$, and the diagonal embeddings $\Delta_k:V_1\rightarrow V_1^k$.
\item{}\textsc{[sections]}
If $W \leq V$ are in $\Ob(\cU)$, then $V/W\in \Ob(\cU)$ and $\Ar(\cU)$ contains the inclusion $\iota:W\rightarrow V$ and quotient $p:V\rightarrow V/W$ maps.
\item{}\textsc{[kernels/images]}
If $f:V_1 \rightarrow V_2$ is in $\Ar(\cU)$, then $\ker f, \im f \in \Ob(\cU)$, and for all $W_1,W_2 \in \Ob(\cU)$,
	\begin{itemize} 
	\item
	if $W_1 \leq \ker f$, the induced map $\overline{f}\colon V_1/W_1 \to V_2$ is in $\Ar(\cU)$;
	\item 
	if $\im f \le W_2 \le V_2$,  the induced map $\check{f}\colon V_1 \to W_2$ is in $\Ar(\cU)$.
	\end{itemize} 
\item{}[$\Z$-\textsc{module structure}]
If $V \in \Ob(\cU)$, then $\Ar(\cU)$ contains the sum  map $\sigma \colon V\times V \to V$ and the multiplication-by-$n$  maps $\mu_n\colon V \to V$.
\end{itemize}
The groups in $\Ob(\cU)$ are called \textbf{modules} of $\cU$ and the homomorphisms in $\Ar(\cU)$ are called \textbf{compatible morphisms}.
If a group $G$ acts on some $V\in \Ob(\cU)$  by {compatible} morphisms, we say that $V$ is a \textbf{$G$-module} in $\cU$.
\end{definition*}

\begin{definition*}\label{d.ModAddDim}
Let $\cU$ be a modular universe.
An \textbf{additive dimension} on $\cU$ is a function $\dim\colon \Ob(\cU) \to \N$ such that for all $f\colon V_1 \to V_2$ in $\Ar(\cU)$, 
\[\dim V = \dim \ker f + \dim \im f.\]
Each $V\in \Ob(\cU)$ is said to be a \textbf{module with an additive dimension}, leaving $\cU$ and $\dim$ implicit from context. Additionally, $V$ is called \textbf{dim-connected}  if for all $W\in \Ob(\cU)$, $W<V$ implies $\dim W < \dim V$.
\end{definition*}

\begin{fact}[{\cite[Connectedness Properties]{CDW23}}]\label{fact.ConnectednessProps}
Let $(\cU, \dim)$ be a modular universe with an additive dimension, and let $V_1,V_2,V\in \Ob(\cU)$.
\begin{enumerate}
\item
For $f\colon V_1 \to V_2$ compatible, if $V_1$ is dim-connected, then so is $\im f$.
\item
If $V_1$ and $V_2$ are dim-connected, then so is $V_1 \times V_2$.
\item
If $V_1, V_2 \leq V$ and $V_1$ and $V_2$ are dim-connected, then so is $V_1+V_2$.
\end{enumerate}
\end{fact} 

\begin{fact}[{\cite[Divisibility Properties, Coprimality Lemma]{CDW23}}]\label{fact.DivisibilityCoprimality}
Let $G$ be a group, and let $V$ be a dim-connected $G$-module with an additive dimension.
\begin{enumerate}
    \item\label{fact.DivisibilityCoprimality.i.Divisibility}
    $V$ is $p$-divisible for $p$ prime if and only if $\dim \{v \in V\mid pv = 0\}=0$. 
    \item\label{fact.DivisibilityCoprimality.i.Coprimality}
    If  $g\in G$ has prime order $p$ and $V$ is $p$-divisible, then  $V = B_g + C_g$ and $\dim (B_g + C_g) = 0$ (sometimes denoted by  $V = B_g \qoplus C_g$).
\end{enumerate}
\end{fact}

\begin{definition*}\label{d.Irreducible}
Let $G$ be a group, and let $V$ be a dim-connected $G$-module with an additive dimension. We say  $V$ is \textbf{dc-irreducible}, if it has no non-trivial, proper, dim-connected $G$-submodule.
\end{definition*}

\begin{fact}[{see \cite[Characteristic Lemma]{CDW23}}]\label{fact.Characteristic}
Let $G$ be a group, and let $V$ be a dc-irreducible $G$-module with an additive dimension. Then either $V$  is divisible or $V$ is an elementary abelian $p$-group for some prime $p$. 
\end{fact}

The main result of~\cite{CDW23} was the classification of the minimal dc-irreducible $\Sym(n)$- and $\Alt(n)$-modules with an additive dimension, assuming $n$ is large enough; here minimal means of minimal dimension for fixed $n$, which is typically $n-2$ or $n-1$. The classification for small $n$ was then picked up in~\cite{CDWY24}. The analysis there determines the minimal dimensions of the modules in question---for all $n$---and also completes the classification of some of the exceptional cases: $\Alt(5)$-modules of dimension $2$ in characteristic $2$, $\Alt(5)$-modules of dimension $3$ in characteristic $5$, and $\Alt(6)$-modules of dimension $3$ in characteristic $3$. One of the main remaining tasks in the classification of the exceptional modules is what began our work on this paper and is addressed by our Corollary~\ref{cor.Alt5Recognition-dim}: to 
classify the minimal faithful $\Alt(5)$-modules in characteristics other than $2$ and $5$ \cite[Remaining Task~1]{CDWY24}.

\begin{proof}[Proof of Corollary~\ref{cor.Alt5Recognition-dim}]
Set $G=\Alt(5)$, and let $V$ be a faithful $G$-module with an additive dimension (as defined above). Further assume that $V$ is dim-connected of dimension $3$ and without $2$-torsion. We verify the hypotheses of Theorem~\ref{thm.Alt5Recognition}. Throughout, we write $C_g$ for $\tr_g(V)$.

Since $V$ is without $2$-torsion, Fact~\ref{fact.DivisibilityCoprimality}\eqref{fact.DivisibilityCoprimality.i.Divisibility} implies  $V$ is $2$-divisible, hence uniquely $2$-divisible. Also, $V$ is dc-irreducible by \cite[Lemmas~3.8 and 3.9]{CDWY24}, so Fact~\ref{fact.Characteristic} shows that $V$ either has prime characteristic or is divisible. 

We next verify the normalizer assumption for $g\in G$ of order $2$. Here,  $K = N_G(\langle g\rangle)$ is a Klein $4$-group, and in the notation of Lemma~\ref{lem.2-divisible-Klein-decomposition}, $V = V_{\alpha_1} \oplus V_{\alpha_2} \oplus V_{\alpha_3} \oplus C_K$ (with $g$ one of the $\alpha_i$). Conjugacy of the $\alpha_i$ ensures that $\dim V_{\alpha_i}$ is independent of $i$, and this dimension must not be $0$ as otherwise $V = C_K$, contradicting faithfulness. Thus  $\dim V_{\alpha_i} = 1$ and, consequently, $C_K$ is trivial (being dim-connected of dimension $0$). We then find---since $V$ is without $2$-torsion---that the definition of the $V_{\alpha_i}$ implies that $V_{\alpha_i} = C_{\alpha_i}$ and hence that $C_{\alpha_i}$ is inverted by $K\setminus \langle \alpha_i\rangle$, as desired.

Finally, take $g\in G$ of order $p\in \{3,5\}$, and assume $V$ is $p$-divisible. Set $N=N_G(\langle g\rangle)$, a dihedral group of order $2p$, and let $\alpha \in N$ be an involution. We show  $\alpha$ inverts $C_g$. 
By Fact~\ref{fact.DivisibilityCoprimality}\eqref{fact.DivisibilityCoprimality.i.Coprimality}, $V=B_g + C_g$ with $\dim(B_g \cap C_g) = 0$. 
Note that $\alpha$ cannot centralize nor invert all of $B_g$ as otherwise $\alpha$ and $g$ would commute on $B_g$ and also on $C_g\le C_V(g)$, hence on $V$. In particular, we must have $\dim B_g \ge 2$ (see \cite[Lemma~3.1]{CDWY24}), so $\dim C_g \le 1$. If $\dim C_g =0$, then $C_g$ is trivial so is certainly inverted by $\alpha$. And if $\dim C_g =1$, then $C_g$ is either centralized or inverted by $\alpha$, but the former is ruled out as it would lead to $C_g = C_\alpha$ and thus $\alpha$ inverting $B_g$, a contradiction.
\end{proof}

\section*{Acknowledgments}
This work benefited greatly from the many helpful comments and suggestions of Adrien Deloro on both early and late versions of this article. The authors are extremely grateful for this as well as for his continued encouragement throughout the project. The authors are also immensely appreciative of the anonymous referee for a highly detailed---and surely time-consuming---report that led to many enhancements, the most notable being a significant streamlining of Proposition~\ref{prop.A5.local}.
The first author would also like to thank her grandmother for being the author's first, greatest cheerleader during life's uncertain trials and all the days in between.

The work of both authors was partially supported by the National Science Foundation under grant No.~DMS-1954127.

\bibliographystyle{alpha}
\bibliography{RotationBib}
\end{document}